\documentclass[12pt]{amsart}

\usepackage{amssymb}
\usepackage{mathrsfs}

\newtheorem{thm}{Theorem}[section]
\newtheorem{lem}[thm]{Lemma}
\newtheorem{prop}[thm]{Proposition}
\newtheorem{claim}[thm]{Claim}
\newtheorem*{shelah_conjecture}{Shelah's Conjecture}
\newtheorem*{01law}{0-1 law for open set mappings}
\newtheorem{subclaim}[thm]{Subclaim}

\newtheorem*{thm*}{Theorem}
\newtheorem{defn}[thm]{Definition}
\newtheorem{question}[thm]{Question}

\theoremstyle{remark}
\newtheorem*{remark}{Remark}

\newcommand\Acal{\mathscr{A}}
\newcommand\Bcal{\mathscr{B}}
\newcommand\Dcal{\mathscr{D}}
\newcommand\Ecal{\mathscr{E}}
\newcommand\Fcal{\mathscr{F}}
\newcommand\Hcal{\mathscr{H}}
\newcommand\Ical{\mathscr{I}}
\newcommand\Ncal{\mathscr{N}}
\newcommand\Pcal{\mathscr{P}}
\newcommand\Qcal{\mathscr{Q}}
\newcommand\Tcal{\mathscr{T}}
\newcommand\Ucal{\mathscr{U}}

\newcommand\Rbb{\mathbb{R}}

\newcommand{\diff}{\operatorname{diff}}
\newcommand{\Ht}{\operatorname{ht}}
\newcommand{\meet}{\wedge}
\newcommand\triord{\triangleleft}
\newcommand{\Th}{{}^{\textrm{th}}}
\newcommand\axiom{\mathrm}
\newcommand\MA{\axiom{MA}}
\newcommand\PFA{\axiom{PFA}}
\newcommand\BPFA{\axiom{BPFA}}
\newcommand\MRP{\axiom{MRP}}
\newcommand\SMRP{\axiom{SMRP}}

\newcommand{\<}{\langle}
\renewcommand{\>}{\rangle}
\newcommand\mand{\textrm{ and }}

\newcommand\fin{\mathrm{fin}}

\title[uncountable linear orders]{A five element basis for the 
uncountable linear orders}

\keywords{Aronszajn, basis,
Countryman,
forcing axiom, linear order,
MRP,
PFA, Shelah's conjecture}

\subjclass[2000]{Primary: 03E35, 03E75, 06A05; Secondary: 03E02 }

\thanks{
This paper is dedicated to Fennel Marie Moore.
I would like to thank J\"{o}rg Brendle for supporting my visit
to Japan (via Grant-in-aid
 for Scientific Research (C)(2)15540120,
Japanese Society for the Promotion of Science)
where I presented the results
of this paper in a series of lectures at Kobe University in
December 2003.
This research was completed while funding for NSF grant DMS--0401893 was
pending and represents the completion of part of the submitted project.
Some revisions and updates were made to the paper after the grant was
funded.
I would like to thank Jean Larson, Paul Larson, Bill Mitchell, and
Boban Veli\v{c}kovi\'{c}
for carefully reading the paper and offering their suggestions
and comments.}
\author{Justin Tatch Moore}

\begin{document}

\begin{abstract}
In this paper I will show that it is relatively consistent with the
usual axioms of mathematics (ZFC) together
with a strong form of the axiom of infinity (the existence of a supercompact
cardinal)
that the class of uncountable linear orders has a five element basis.
In fact such a basis follows from the Proper Forcing Axiom, a strong form
of the Baire Category Theorem.
The elements are $X$, $\omega_1$, $\omega_1^*$, $C$, $C^*$ where $X$ is
any suborder of the reals of cardinality $\aleph_1$ and $C$ is any Countryman line.
This confirms a longstanding conjecture of Shelah.
\end{abstract}

\maketitle

\section{Introduction}

Our focus in this paper will be to show that the Proper Forcing
Axiom ($\PFA$) implies
that any uncountable linear order
must contain an isomorphic copy of one of the following five orders
$X$, $\omega_1$, $\omega_1^*$, $C$, and $C^*$.
Here $X$ is any fixed set of reals of cardinality $\aleph_1$ and $C$ is any fixed Countryman line.
Such a list is called a \emph{basis}.

The simplest example of an uncountable linear order is $\Rbb$, the real line.
This object --- familiar to every mathematician ---
serves as the prototype for the class of linear orders and
as the canonical example of an uncountable set.
Early on in modern set theory Baumgartner proved the following deep result
which suggested that it might be possible to prove more general
classification results for uncountable linear orders.
\begin{thm}\cite{reals_isomorphic}
($\PFA$)
If two sets of reals are $\aleph_1$-dense,\footnote{
I.e. every interval meets them at a set of cardinality $\aleph_1$.}
then they are isomorphic.
In particular if $X$ is a set of reals of cardinality $\aleph_1$, then
$X$ serves as a single element basis for the class of uncountable
separable linear orders.
\end{thm}
$\PFA$ is a strengthening
of Baire's Category Theorem and is independent of the usual axioms of set
theory.
Its use in infinite combinatorics can be likened to Erd\"os's
\emph{probabilistic method} (see \cite{probabilistic_method}) from
finite combinatorics.
The main differences are that the notion of a probability space
is replaced by the more abstract notion of a proper forcing and
the assertion ``If an object can be chosen with positive probability, then
it exists'' requires an axiomatic assumption.
Frequently --- as in Baumgartner's result above --- this axiom can
be used to find morphisms between certain structures
or to make other combinatorial reductions
(see \cite{club_isomorphic},
\cite{reals_isomorphic},
\cite{class_trans},
\cite{OCA_aut}).

Some additional assumption is necessary in Baumgartner's result
because of the following classical construction of Sierpi\'{n}ski.
\begin{thm}\cite{many_iso_types}
There is a set of reals $X$ of cardinality continuum such that
if $f \subseteq X^2$ is a continuous injective function, then
$f$ differs from the identity function on a set
of cardinality less than continuum. 
\end{thm}
From this it is routine to prove that under
the Continuum Hypothesis there is no basis for the uncountable
separable linear orders of cardinality less than $|\Pcal(\Rbb)|$.
This gives a complete contrast to the conclusion of Baumgartner's result.

The simplest example of a linear order which is separable only in the
trivial instances is a well order.
The uncountable well orders have a canonical
minimal representative, the ordinal $\omega_1$.\footnote{The canonical
representation of well orders mentioned here is due to von Neumann.}
Similarly, the converse $\omega_1^*$ of $\omega_1$ obtained by reversing the
order relation forms
a single element basis for all of the uncountable converse well orders.

Those uncountable linear orders which do not contain uncountable
separable suborders or copies of $\omega_1$ or $\omega_1^*$ are
called \emph{Aronszajn lines}.\footnote{Or \emph{Specker types.}}
They are classical objects considered long ago by Aronszajn who first
proved their existence.
Some time later Countryman made a
brief but important contribution to the subject by
asking whether there is an uncountable linear order $C$ whose square
is the union of countably many chains.\footnote{Here
\emph{chain} refers to the coordinate-wise partial order on $C^2$.}
Such an order is necessarily Aronszajn.
Furthermore, it is easily seen that no uncountable linear order can
embed into both a Countryman line and its converse.
Shelah proved that such orders exist in ZFC \cite{Countryman:Shelah}
and made the following conjecture:\footnote{Before this point is
was an open problem whether the uncountable
linear orders had a four element basis.
Also, Shelah simply conjectured the consistency of such a basis.
This was at least in part because the language of proper forcing
and $\PFA$ was not around at the time.
Still, it is very reasonable to assume that is is how the conjecture
would have been phrased had the language been available and certainly
this is how the conjecture was viewed by the end of the 1980's.}
\begin{shelah_conjecture}\cite{Countryman:Shelah} ($\PFA$)
The orders
$X$, $\omega_1$, $\omega_1^*$, $C$ and $C^*$ form
a five element basis for the uncountable linear orders any time
$X$ is a set of reals of cardinality $\aleph_1$ and $C$ is
a Countryman line.
\end{shelah_conjecture}
\noindent
Notice that by our observations
such a basis is necessarily minimal.

This problem was exposited, along with some other basis problems for
uncountable structures, in Todor\v{c}evi\'{c}'s address
to the 1998 International Congress of Mathematicians \cite{basis_problems}.
It also appears as Question 5.1 in Shelah's problem list
\cite{problems:Shelah}.
In this paper I will prove Shelah's conjecture.
In doing so, I will introduce some new methods for applying
$\PFA$ which may be relevant to solving other problems.

\section{Background}

This paper should be readily accessible to
anyone who is well versed in set theory and
the major developments in the field in the 70's and 80's.
The reader is assumed to have proficiency in the areas of
Aronszajn tree combinatorics, forcing axioms,
the combinatorics of $[X]^{\aleph_0}$, and Skolem hull arguments.
Jech's \cite{set_theory:Jech} and Kunen's \cite{set_theory:Kunen} serve as
good references on general set theory.
They both contain some basic information on Aronszajn trees;
further reading on Aronszajn trees can be found in \cite{trees:Todorcevic} and
\cite{cseq}.The reader is referred to
\cite{set_theory:Bekkali}, \cite{Countryman:Shelah}, \cite{l-map},
\cite{partitioning_ordinals},
or \cite{cseq} for information
on Countryman lines.
It should be noted, however, that knowledge of the method of minimal walks will
not be required.
The set theoretic assumption we will be working with is
the Proper Forcing Axiom.
Both \cite{class_trans}
and the section on PFA in \cite{partition_problems}
serve as good concise references on the subject for our purposes.
A more elaborate account of proper forcing can be found in Shelah's
\cite{proper_forcing}.
See \cite{MRP} for information on the Mapping Reflection Principle.
For basic forcing technology, the reader is referred to
\cite{multiple_forcing} and \cite{set_theory:Kunen}.
Part III of Jech's \cite{multiple_forcing} gives a good
exposition on the combinatorics of
$[X]^{\aleph_0}$, the corresponding closed unbounded (or \emph{club}) filter,
and related topics.

The notation in this paper is mostly standard.
If $X$ is an uncountable set, then $[X]^{\aleph_0}$ will
be used to denote the collection of all countable
subsets of $X$.
All ordinals are von Neumann ordinals --- they are the set
of their predecessors under the $\in$ relation.
The collections $H(\theta)$ for regular cardinals $\theta$ consist of
those sets of hereditary cardinality less than $\theta$.
Hence $H({2^{\theta}}^+)$ contains $\Pcal(H(\theta^+))$ as a subset
and $H(\theta)$ as an element.
Often when I refer to $H(\theta)$ in this paper I will really be
referring to the structure
$(H(\theta),\in,\triord)$ where $\triord$ is some fixed well ordering of
$H(\theta)$ which
can be used to generate the Skolem functions.
 
\section{The axioms}

The working assumption in this paper will be
the Proper Forcing Axiom introduced by
Shelah and proved relatively consistent from a supercompact cardinal.
We will often appeal to the bounded form of this axiom
isolated by Goldstern and Shelah \cite{BPFA}.
We will use an equivalent formulation due to Bagaria \cite{gen_abs_FA}:
\begin{description}

\item[$\BPFA$]
If $\phi$ is a formula in language of $H(\aleph_2)$ with only
bounded quantifiers
and there is a proper partial order which forces $\exists X \phi(X)$, then
$H(\aleph_2)$ already satisfies $\exists X \phi(X)$.

\end{description}
At a crucial point in the proof we will also employ the
Mapping Reflection Principle introduced recently in \cite{MRP}.
In order to state it we will need the following definitions.
\begin{defn}
If $X$ is an uncountable set, then there is a natural topology
--- the Ellentuck topology ---
on $[X]^{\aleph_0}$ defined by declaring
$$[x,N] = \{Y \in [X]^{\aleph_0}:x \subseteq Y \subseteq N\}$$
to be open whenever $N$ is in
$[X]^{\aleph_0}$ and $x$ is a finite subset of $N$.
\end{defn}

This topology is regular and 0-dimensional.
Moreover, the closed and cofinal sets generate the club filter
on $[X]^{\aleph_0}$.

\begin{defn}
If $M$ is an elementary submodel of some $H(\theta)$ and $X$
is in $M$, then we say a subset $\Sigma \subseteq [X]^{\aleph_0}$
is \emph{$M$-stationary} if whenever $E \subseteq [X]^{\aleph_0}$
is a club in $M$, the intersection
$\Sigma \cap E \cap M$ is non-empty.
\end{defn}
\begin{defn}
If $\Sigma$ is a set mapping defined on a set of countable
elementary submodels of some $H(\theta)$ and there is an $X$ such that
$\Sigma(M) \subseteq [X]^{\aleph_0}$ is open and $M$-stationary
for all $M$, then we say $\Sigma$ is an
\emph{open stationary set mapping}.
\end{defn}
The Mapping Reflection Principle is the following statement:
\begin{description}

\item[$\MRP$]
If $\Sigma$ is an open stationary set mapping defined on a club of
models, then there is a continuous $\in$-chain
$\<N_\xi:\xi < \omega_1\>$ in the domain of $\Sigma$ such that for
every $\nu > 0$ there is a $\nu_0 < \nu$ such that
$N_\xi \cap X$ is in $\Sigma(N_\nu)$ whenever
$\nu_0 < \xi < \nu$.

\end{description}
The sequence $\<N_\xi:\xi < \omega_1\>$ postulated by this axiom
will be called a
\emph{reflecting sequence} for the set mapping $\Sigma$.

\section{A combinatorial reduction}

Rather than prove Shelah's basis conjecture directly, I will appeal
to an observation of Abraham and Shelah in \cite{club_isomorphic}.
\begin{thm} \label{reduction_thm} \cite{club_isomorphic} ($\BPFA$)
The following are equivalent:
\begin{enumerate}

\item The uncountable linear orders have a five element basis.

\item There is an Aronszajn tree $T$ such that for every
$K \subseteq T$ there is an uncountable antichain $X \subseteq T$ such
that $\meet(X)$\footnote{This will be defined momentarily.}
is either contained in or disjoint from $K$.

\end{enumerate}
\end{thm}
A detailed proof of this theorem can be found in the last section of
\cite{l-map}.
I will sketch the proof for completeness.

The implication \emph{(1) implies (2)} does not require $\BPFA$ and in fact
(1) implies that the conclusion of (2) holds for an arbitrary
Aronszajn tree $T$.
To see why it is true, suppose that $(T,\leq)$ is an Aronszajn tree
equipped with a lexicographical order
and suppose that $K \subseteq T$ witnesses a failure of (2).
If $(T,\leq)$ doesn't contains a Countryman suborder,
then (1) must fail.
So without loss of generality, we may assume that $(T,\leq)$ is
Countryman.

Define $s \leq' t$ iff $s \meet t$ is in $K$ and $s \leq t$ or
$s \meet t$ is not in $K$ and $t \leq s$.
It is sufficient to check that neither $(T,\leq)$ nor its converse
$(T,\geq)$ embeds an uncountable suborder of $(T,\leq')$.
This is accomplished with two observations.
First, since $(T,\leq)$ and its converse are Countryman, any
such embedding can be assumed to be the identity map.
Second, if $\leq$ and $\leq'$ agree
on $X \subseteq T$, then $\meet(X) \subseteq K$;
disagreement on $X$ results in $\meet(X) \cap K = \emptyset$.

For the implication \emph{(2) implies (1)} we first observe that,
by Baumgartner's result mentioned above,
it suffices to show that the Aronszajn lines have a two element basis.
Fix a Countryman line $C$ which is a lexicographical
order $\leq$
on an Aronszajn tree $T$.
The club isomorphism of Aronszajn trees
under $\BPFA$ \cite{club_isomorphic} together with
some further appeal to $\MA_{\aleph_1}$ implies that any Aronszajn line
contains a suborder isomorphic to some $(X,\leq')$ where
$X \subseteq T$ is uncountable and binary and $\leq'$ is a --- possibly
different --- lexicographical order on $T$.
Statement (2) is used to compare $\leq$ and $\leq'$ and find an
uncountable $Y \subseteq X$ on which they always agree or always disagree.
Applying $\MA_{\aleph_1}$, $C$ embeds into all its uncountable
suborders, thus finishing the proof.

\section{The proof of the main result}

In this section we will prove the basis conjecture of Shelah by
proving the following result and appealing
to Theorem \ref{reduction_thm}.

\begin{thm} ($\PFA$)
There is an Aronszajn tree $T$ such that 
if $K \subseteq T$, then
there is an uncountable antichain $X \subseteq T$ such that
$\meet(X)$ is either contained in or disjoint from $K$.
\end{thm}

The proof will be given as a series of lemmas.
In each case, I will state any set theoretic hypothesis needed to
prove a lemma.
This is not so much to split hairs but because I feel that it will help
the reader better understand the proof.

For the duration of the proof, we will let $T$ be a fixed Aronszajn tree
which is contained in the complete binary tree, coherent,
closed under finite changes, and special.\footnote{The tree
$T(\rho_3)$ of \cite{partitioning_ordinals} is such an example.}
It will be convenient to first make some definitions and fix some
notation.

\begin{defn} If $s$ and $t$ are two elements of $T$, then $\diff(s,t)$
is the set of all $\xi$ such that $s(\xi) \ne t(\xi)$.
If $F \subseteq T$, then $\diff(F)$ is the union of all $\diff(s,t)$ such
that $s$ and $t$ are in $F$.\footnote{
Coherence is just the assertion that $\diff(s,t)$ is a finite set for all
$s,t$ in $T$.}
\end{defn}

\begin{defn}
If $X$ is a subset of $T$ and $\delta < \omega_1$, then
$X \restriction \delta$ is the set of all $t \restriction \delta$ such
that $t$ is in $X$.
Here $t \restriction \delta$ is just functional restriction.
\end{defn}

\begin{defn}
If $s$ and $t$ are in $T$, then $\Delta(s,t)$ is the least
element of $\diff(s,t)$.
If $s$ and $t$ are comparable, we leave
$\Delta(s,t)$ undefined.\footnote{This is somewhat non-standard but
it will simplify the notation at some points.
For example, in the definition of
$\Delta(Z,t)$ we only collect those
values where $\Delta$ is defined.}
If $Z \subseteq T$ and $t$ is in $T$, then
$\Delta(Z,t) = \{\Delta(s,t): s \in Z\}$.
\end{defn}

\begin{defn}
If $X$ is a finite subset of $T$, then $X(j)$ will denote the $j\Th$ least
element of $X$ in the lexicographical order inherited from $T$.
\end{defn}

\begin{defn}
If $s,t$ are incomparable in $T$, then the meet of $s$ and $t$ --- denoted
$s \meet t$ --- is the restriction
$s \restriction \Delta(s,t) = t \restriction \Delta(s,t)$.
If $X$ is a subset of $T$, then
$\meet(X) = \{s \meet t : s,t \in X\}$.\footnote{
The domain of $\meet$ is the same as the domain of $\Delta$; the
set of all incomparable pairs of elements of $T$.}
\end{defn}

The following definition provides a useful means of measuring
subsets of an elementary submodel's intersection with $\omega_1$.

\begin{defn}
If $P$ is a countable elementary submodel of $H(\aleph_2)$ containing
$T$ as an element, define $\Ical_P(T)$ to be the collection of all
$I \subseteq \omega_1$ such that for some uncountable $Z \subseteq T$
in $P$ and some $t$ of height $P \cap \omega_1$
which is in the downward closure of $Z$, the set
$\Delta(Z,t)$ is disjoint from $I$.
\end{defn}

The following propositions are routine to verify using the
coherence of $T$ and its closure under finite changes
(compare to the proof that $\Ucal(T)$ is a filter in
\cite{l-map} or \cite{cseq}).

\begin{prop}\label{I_P(T)_is_ideal}
If $I$ is in $\Ical_P(T)$ and $t$ is in $T$ with height $P \cap \omega_1$,
then there is a $Z \subseteq T$ in $P$ such that
$t$ is in the downward closure\footnote{
The \emph{downward closure} of $Z$ is the collection of all
$s$ such that $s \leq s^*$ for some $s^*$ in $Z$}
of $Z$ and $\Delta(Z,t)$ is disjoint from $I$.
\end{prop}

\begin{prop}\label{I_P(T)_refine}
If $I$ is in $\Ical_P(T)$, $Z_0$ is a subset of $T$ in $P$ and
$t$ is an element of the downward closure of $Z_0$ of height
$P \cap \omega_1$, then there is a $Z \subseteq Z_0$ in $P$ which
also contains $t$ in its downward closure and satisfies
$\Delta(Z,t) \cap I$ is empty.
\end{prop}

\begin{prop}
$\Ical_P(T)$ is a proper ideal on $\omega_1$ which contains
$I \subseteq \omega_1$ whenever $I \cap P$ is bounded in
$\omega_1 \cap P$.
\end{prop}

\begin{prop}
Suppose $P$ is a countable elementary submodel of $H(\aleph_2)$ such that
$Z \subseteq T$ is an element of
$P$, and there is a $t \in T$ of height $P \cap \omega_1$
in the downward closure of $Z$. 
Then $Z$ is uncountable.
\end{prop}

Let $K \subseteq T$ be given.
The following definitions will be central to the proof.
The first is the na\"\i ve approach to
forcing an uncountable $X$ such that $\meet(X)$ is contained in $K$.

\begin{defn}
$\Hcal(K)$ is the collection of all finite $X \subseteq T$ such
that $\meet(X)$ is contained in $K$.\footnote{
A collection of finite sets such as this becomes a forcing
notion when given the order of reverse inclusion ($q \leq p$ means
that $q$ is stronger than $p$).
A collection of ordered pairs of finite sets becomes a forcing by
coordinate-wise reverse inclusion.}
\end{defn}

The second is the notion of rejection which will be central
in the analysis of $\Hcal(K)$.
For convenience we will let $\Ecal$ denote the collection of all
clubs $E \subseteq [H(\aleph_2)]^{\aleph_0}$ which consist of
elementary submodels which contain $T$ and $K$ as elements.
Let $E_0$ denote the element of $\Ecal$ which consists
of all such submodels.
\begin{defn}
If $X$ is a finite subset of $T$, then
let $K(X)$ denote the set of all $\gamma < \omega_1$ such that
for all $t$ in $X$, if $\gamma$ is less than the height of $t$, then
$t \restriction \gamma$ is in $K$.
\end{defn}

\begin{defn}
If $P$ is in $E_0$ and $X$ is a finite subset of $T$, then
we say that \emph{$P$ rejects $X$} if
$K(X)$ is in $\Ical_P(T)$.
\end{defn}

The following trivial observations about $P$ in $E_0$
and finite $X \subseteq T$ are useful and will be used tacitly
at times in the proofs which follow.

\begin{prop}
If $P$ does not reject $X$, then it
does not reject any of its restrictions.
\end{prop}

\begin{prop} \label{rejection_prop}
$P$ rejects $X$ iff it rejects $X \restriction (P \cap \omega_1)$
iff it rejects $X \setminus P$.
\end{prop}

\begin{prop}
If $X$ is in $P$, then $P$ does not reject $X$.
\end{prop}

The forcing notion $\partial (K)$ which we are about to define
seeks to add a subset of $T$ in which
rejection is rarely encountered.\footnote{
The symbol $\partial$ is being used here
because there is a connection to the notion of a Cantor-Bendixon derivative.
In a certain sense we are removing the parts of the
partial order $\Hcal(K)$ which
are causing it to be improper.}

\begin{defn}
$\partial (K)$ consists of all pairs $p = (X_p,\Ncal_p)$
such that:
\begin{enumerate}

\item $\Ncal_p$ is a finite $\in$-chain of countable elementary
submodels of $H({2^{\aleph_1}}^+)$ each of which
contain $T$, $K$, and $E_0$ as members.

\item \label{stationary_positives}
$X_p \subseteq T$ is a finite set and
if $N$ is in $\Ncal_p$, then there is an $E$ in $\Ecal \cap N$ such
that $X_p$ is not rejected by any element of $E \cap N$.

\end{enumerate}
We will also be interested in the suborder
$$\partial\Hcal(K) = \{p \in \partial (K): X_p \in \Hcal(K)\}$$
which seems to be the correct modification of $\Hcal(K)$ from
the point of view of forcing the conclusion of the main theorem.
\end{defn}
In order to aid in the presentation of the lemmas,
I will make the following definition.
\begin{defn}
$\partial (K)$ is \emph{canonically proper} if
whenever $M$ is a countable elementary submodel of
$H\big(|2^{\partial(K)}|^+\big)$ and $\partial(K)$ is in $M$,
any condition $p$ which satisfies $M \cap H({2^{\aleph_1}}^+)$ is in
$\Ncal_p$ is $(M,\partial (K))$-generic.
An analogous definition is made for $\partial\Hcal(K)$.
\end{defn}
We will eventually prove that, assuming the Proper Forcing Axiom,
$\partial\Hcal(K)$ is canonically proper.
The following lemma shows that this is sufficient to finish the argument.

\begin{lem} ($\BPFA$)
If $\partial \Hcal(K)$ is canonically proper, then
there is an uncountable $X \subseteq T$ such that
$\meet(X)$ is either contained in $K$ or disjoint from $K$.
\end{lem}

\begin{remark}
This conclusion is sufficient since the properties of $T$ imply
that $X$ contains an uncountable antichain.
\end{remark}

\begin{proof}
Let $M$ be an elementary submodel of $H\big(|2^{\partial \Hcal(K)}|^+\big)$
containing $\partial \Hcal(K)$ as an element.
Let $t$ be an element of $T$ of height $M \cap \omega_1$.
If $$p = \big(\{t\},\{M \cap H({2^{\aleph_1}}^+)\}\big)$$
is a condition in $\partial \Hcal(K)$, then it is
$(M,\partial \Hcal(K))$-generic by assumption.
Consequently $p$ forces that the interpretation of
$$\dot X = \{s \in T: \exists q \in \dot G (s \in X_q)\}$$
is uncountable.
Since $\dot X$ will then be forced to have the property
that $\meet(\dot X) \subseteq \check K$, we can apply $\BPFA$
to find such an $X$ in $V$.

Now suppose that $p$ is not a condition.
It follows that there is a countable elementary submodel $P$ of $H(\aleph_2)$
in $M$ such that $T$ is in $P$ and $K(\{t\})$
is in $\Ical_P(T)$.
Therefore there is a $Z \subseteq T$ in $P$ such that
$t \restriction (P \cap \omega_1)$ is in the downward closure
of $Z$ and for all $s$ in $Z$,
$s \meet t$ is not in $K$. 
Let $Y$ consist of all those $w$ in $\meet(Z)$ such that if $u,v$ are
incomparable elements of $Z$ and $u \meet v \leq w$, then
$u \meet v$ is not in $K$.
Notice that $Y$ is an element of $P$.
$Y$ is uncountable since it contains $s \meet t$ for every $s$ in
$P \cap Z$ which is incomparable with $t$ and the heights of elements
of this set is easily seen to be unbounded in $P \cap \omega_1$.
We are therefore finished once we see that $\meet(Y)$ is disjoint from $K$.
To this end, suppose that $w_0$ and $w_1$ are incomparable elements of $Y$.
Let $u_0,u_1,v_0,v_1$ be elements of $Z$ such that
$u_i$ and $v_i$ are incomparable and $w_i = u_i \meet v_i$.
Since $w_0$ and $w_1$ are incomparable,
$$u_0\big(\Delta(w_0,w_1)\big) = w_0\big(\Delta(w_0,w_1)\big) \ne
w_1\big(\Delta(w_0,w_1)\big) = v_1\big(\Delta(w_0,w_1)\big).$$
It follows that $u_0 \meet v_1 = w_0 \meet w_1$.
Since $w_0$ extends $u_0 \meet v_1$ and is in $Y$, it must
be that $u_0 \meet v_1$ is not in $K$.
Hence $w_0 \meet w_1$ is not in $K$.
This completes the proof that $\meet(Y)$ is disjoint from $K$.
\end{proof}

The following lemma is the reason for our definition of
rejection.
It will be used at crucial points in the argument.

\begin{lem}
\label{rel_c.c.c.}
Suppose that $E$ is in $\Ecal$
and $\<X_\xi:\xi < \omega_1\>$ is a sequence of disjoint $n$-element
subsets of $T$
so that no element of $E$ rejects any $X_\xi$ for $\xi < \omega_1$.
Then there are $\xi \ne \eta < \omega_1$ such that
$X_\xi(j) \meet X_\eta(j)$ is in $K$ for all $j < n$.
\end{lem}

\begin{proof}
By the pressing down lemma we can find a $\zeta < \omega_1$ and a stationary
set $\Xi \subseteq \omega_1$ such that:
\begin{enumerate}

\item For all $\xi$ in $\Xi$, $X_\xi$ contains only elements
of height at least $\xi$.

\item 
$X_\xi(j) \restriction \zeta = X_\eta(j) \restriction \zeta$
for all $j < n$ and $\xi,\eta \in \Xi$.

\item For all $\xi$ in $\Xi$ the
set $\diff(X_\xi \restriction \xi)$ is contained in $\zeta$.

\end{enumerate}
Now let $P$ be an element of $E$ which contains
$\<X_\xi:\xi \in \Xi\>$.
Let $\eta$ be an element of $\Xi$ outside of $P$ and pick
a $\xi$ in $\Xi \cap P$ such that
$X_\xi(0) \restriction \xi$ and $X_\eta(0) \restriction \eta$
are incomparable and
for all $j < n$
$$X_\eta(j) \restriction \Delta\big(X_\eta(0),X_\xi(0)\big)$$
is in $K$.
This is possible since otherwise
$Z = \{X_\xi(0) \restriction \xi: \xi \in \Xi \}$ and
$t = X_\eta(0) \restriction (P \cap \omega_1)$ would
witness
$K(X_\eta)$
is in $\Ical_P(T)$ and therefore that $P$ rejects $X_\eta$.

Notice that if $j < n$, then
$$\Delta\big(X_\eta(j),X_\xi(j)\big) = \Delta\big(X_\eta(0),X_\xi(0)\big)$$
since
$$\diff(X_\xi \restriction \xi) \cup \diff(X_\eta \restriction \eta)
\subseteq \zeta,$$
$$X_\eta(j) \restriction \zeta = X_\xi(j) \restriction \zeta.$$
Hence the meets
$$X_\xi(j) \meet X_\eta(j) = X_\eta(j) \restriction \Delta\big(X_\xi(0),X_\eta(0)\big)$$
are in $K$ for all $j < n$.
\end{proof}

The next lemma draws the connection between $\partial \Hcal(K)$ and the forcing
$\partial (K)$.
We will then spend the remainder of the paper analyzing $\partial(K)$.

\begin{lem} ($\BPFA$)
If $\partial(K)$ is canonically proper, so is
$\partial\Hcal(K)$.
\end{lem}

\begin{proof}
We will show that otherwise the forcing $\partial (K)$ introduces
a counterexample to Lemma \ref{rel_c.c.c.} which would then
exist in $V$ by an application of $\BPFA$.
Let $M$ be a countable elementary submodel
of $H\big(|2^{\partial (K)}|^+\big)$ which contains
$K$ as an element and let $r \in \partial \Hcal(K)$ be such
that $M \cap H({2^{\aleph_1}}^+)$ is in $\Ncal_r$ and
yet $r$ is not $(M,\partial \Hcal(K))$-generic.
By extending $r$ if necessary, we may assume that
there is a dense open set $\Dcal \subseteq \partial \Hcal(K)$ in $M$ 
which contains $r$ such that if $q$ is in $\Dcal \cap M$, then
$q$ is $\partial\Hcal(K)$-incompatible with $r$.

Let $E \in \Ecal \cap M$ be such that no element of $E \cap M$ rejects
$X_r$ and let $E'$ be the
elements of $E$ which are the union of their intersection with $E$.
Put $Y_r = (X_r \setminus M) \restriction (M \cap \omega_1)$.
\begin{claim} \label{step_up}
No element of $E'$ rejects $Y_r$.
\end{claim}
\begin{proof}
Let $P$ be an element of $E'$.
We need to verify that
$K(Y_r)$
is not in $\Ical_P(T)$.
If $P \cap \omega_1$ is greater than $M \cap \omega_1$, then
$Y_r \subseteq P$ and this is trivial.
Now suppose that $Z \subseteq T$ is in $P$ and $t$ is an
element of $T$ of height $P \cap \omega_1$ which is in the
downward closure of $Z$.
Let $P_0$ be an element of $E \cap P$ which contains $Z$ as a member.
Such a $P_0$ will satisfy
$$P_0 \cap \omega_1 < P \cap \omega_1 \leq M \cap \omega_1.$$
Let $\nu = P_0 \cap \omega_1$.
If $\Delta\big(Z,t \restriction (P \cap \omega_1)\big)$
is disjoint from $K(Y_r)$, then it witnesses
$K(Y_r \restriction \nu)$ is in $\Ical_{P_0}(T)$.
But then we could use the elementarity of $M$ to find such a $P_0$ in
$M \cap E$,
which is contrary to our choice of $E$.
Hence no element of $E'$ rejects $Y_r$.
\end{proof}
Let $\zeta \in M \cap \omega_1$ be an upper bound for
$\diff(Y_r)$ and let $n = |Y_r|$.
If $j < n$, let $A_j \subseteq T$ be an antichain in $M$ which
contains $Y_r(j)$.
Put $\Dcal_*$ to be the collection of all $q$ in $\Dcal$ such that
\begin{enumerate}

\item $X_r \cap M = X_q \cap N(q)$ where $N(q)$ is the least
element of $\Ncal_q$ which is not in $\Ncal_r \cap M$.

\item
$Y_q \restriction \zeta = Y_r \restriction \zeta$ where
$Y_q = \big(X_q \setminus N(q)\big) \restriction \big(N(q) \cap \omega_1\big)$.

\item \label{incompatible}
No element of $E'$ rejects $Y_q$.

\item
$Y_q(j)$ is in $A_j$ whenever $j < n$.

\end{enumerate}
Note that $\Dcal_*$ is in $M$.

Let $G$ be a $\partial (K)$-generic filter which contains $r$.
Notice that $r$ is $(M,\partial (K))$-generic.
Working in $V[G]$, let $\Fcal$ be the collection of all
$Y_q$ where $q$ is in $\Dcal_* \cap G$.
Now $M[G \cap M]$ is an elementary submodel of
$H\big(|2^{\partial (K)}|^+\big)[G]$\footnote{By Theorem 2.11 of
\cite{proper_forcing}.}
which contains $\Fcal$ as an element but not as a subset
(since $Y_r$ is in $\Fcal$).
Therefore $\Fcal$ is uncountable.
Notice that every element of $\Fcal$ has the property that
it is in $\Hcal(K)$ but that
for every countable $\Fcal_0 \subseteq \Fcal$ there is a $Y_q$ in
$\Fcal \setminus \Fcal_0$ such that $Y_q \cup Y_{q_0}$ is not in $\Hcal(K)$ 
for any $Y_{q_0}$ in $\Fcal_0$.
This follows from the elementarity of $M[G \cap M]$ and
from the fact that $Y_r \cup Y_q$ is not in $\Hcal(K)$ for any $Y_q$ in
$\Fcal \cap M[G \cap M]$.
Now it is possible to build an uncountable sequence
$\<X_\xi : \xi \in \Xi\>$
of elements of $\Fcal$ such that:
\begin{enumerate}

\item
$X_\xi$ has size $n$ for all $\xi \in \Xi$ and is a subset of the
$\xi\Th$ level of $T$.

\item \label{H(K)_incompatible:cond}
$X_\xi \cup X_\eta$ is not in $\Hcal(K)$ whenever $\xi \ne \eta \in \Xi$.

\item \label{same_restriction:cond}
There is a $\zeta < \omega_1$ such that
$X_\xi \restriction \zeta = X_\eta \restriction \zeta$ has size $n$
for all $\xi, \eta < \omega_1$.

\end{enumerate}
It follows from item \ref{H(K)_incompatible:cond}
that if $\xi < \eta < \omega_1$, then there are
$j,j' < n$ such that
$X_\xi(j) \meet X_\eta(j')$ is not in $K$.
By item \ref{same_restriction:cond},
it must be the case that $j = j'$ since
this condition ensures that
$$X_\xi(j) \meet X_\eta(j') = X_\xi(j) \meet X_\xi(j')$$
whenever $j \ne j' < n$ and
hence this meet would be in $K$ by virtue of $X_\xi$ being in $\Hcal(K)$.
Applying $\BPFA$ we get a sequence of sets 
satisfying 1--3 in $V$ and therefore
a contradiction to Lemma \ref{rel_c.c.c.} since no elements of $\Fcal$
are rejected by any member of $E'$.
Hence $\partial \Hcal(K)$ must also be canonically proper.
\end{proof}

Next we have a typical ``models as side conditions'' lemma.

\begin{lem}\label{A_B:lem}
If $\partial (K)$ is not canonically proper, then there are disjoint sets
$\Acal$,
$\Bcal$ and a function
$Y:\Acal \cup \Bcal \to [T]^{<{\aleph_0}}$ such that
\begin{enumerate}

\item $\Acal \subseteq [H({2^{\aleph_1}}^+)]^{{\aleph_0}}$,
$\{N \cap H(\aleph_2): N \in \Acal\}$ is stationary, and every element
of $\Acal$ is an intersection of an elementary submodel of
$H({2^{2^{\aleph_1}}}^+)$ with $H({2^{\aleph_1}}^+)$.

\item $\Bcal \subseteq [H({2^{2^{\aleph_1}}}^+)]^{\aleph_0}$ is stationary.

\item If $M$ is in $\Acal \cup \Bcal$, then
$\big(Y(M),\{M \cap H({2^{\aleph_1}}^+)\}\big)$ is in $\partial (K)$.
 
\item For every $M$ in $\Bcal$ and $N$ in $\Acal \cap M$,
$\big(Y(N) \cup Y(M),\{N\}\big)$ is not a condition in $\partial (K)$.

\end{enumerate}
\end{lem}

\begin{proof}
Let $M$ be a countable elementary submodel of
$H\big(|{2^{\partial (K)}}|^+\big)$ and
$r$ in $\partial (K) \cap M$ be a condition which is not
$(M,\partial (K))$-generic
such that $M \cap H({2^{\aleph_1}}^+)$ is
in $\Ncal_r$.
By extending $r$ if necessary, we can find dense open
$\Dcal \subseteq \partial (K)$ in $M$ which contains $r$ such
that no element of $\Dcal \cap M$ is compatible with $r$.
Furthermore we may assume that if
$q$ is in $\Dcal$, $N$ is in $\Ncal_q$, and $t$ is in $X_q$,
then $t \restriction (N \cap \omega_1)$ is also in $X_q$.

Define $r_0 = (X_r \cap M,\Ncal_r \cap M)$.
If $q$ is in $\Dcal$, let $N(q)$ be the $\in$-least element of
$\Ncal_q \setminus \Ncal_{r_0}$.
Let $k = |\Ncal_r \setminus \Ncal_{r_0}|$
and $\zeta$ be the maximum of all ordinals of the forms
$\Ht(s) + 1$ for $s \in X_{r_0}$
and $N \cap \omega_1$ for $N \in \Ncal_{r_0}$. 
Let $\Tcal_k$ be the set of all $q$ in $\Dcal$ such that:
\begin{enumerate}

\item $X_q \cap N(q) = X_{r_0}$ and $\Ncal_q \cap N(q) = \Ncal_{r_0}$.

\item For all $N$ in $\Ncal_q$,
$N$ is an intersection of an elementary submodel of
$H({2^{2^{\aleph_1}}}^+)$ with $H({2^{\aleph_1}}^+)$.

\item $X_q \restriction \zeta = X_r \restriction \zeta$.

\item $|X_q| = |X_r| = m$ and $|\Ncal_q \setminus \Ncal_{r_0}| = k$.
\end{enumerate}
Let $N_i(q)$ denote the $i\Th$ $\in$-least element of $\Ncal_q \setminus N(q)$
and define $\Tcal_i$ recursively for $i \leq k$.
Given $\Tcal_{i+1}$, define $\Tcal_i$ to be the collection of
all $q$ such that
$$\{N_{i+1}(q^*) \cap H(\aleph_2): q^* \in \Tcal_{i+1} \mand
q = q^* \restriction N_{i+1}(q^*)\}$$
is stationary where
$$q^* \restriction N_{i+1}(q^*)=
\big(X_{q^*} \cap N_{i+1}(q^*),\Ncal_{q^*} \cap N_{i+1}(q^*)\big).$$
Let $\Tcal$ be the collection of all $q$ in
${\displaystyle \bigcup_{i \leq k} \Tcal_i}$
such that if $q$ is in $\Tcal_i$, then $q \restriction N_{i'+1}(q)$
is in $\Tcal_{i'}$ for all $i' < i$.

\begin{claim} 
$r$ is in $\Tcal$.
\end{claim}

\begin{proof}
If $q$ is in $\partial (K)$, define
$$\tilde q = \big(X_q,\{N \cap H(\aleph_2):N \in \Ncal_q\}\big).$$
While elements of $\Ncal_r \setminus M$ need not contain
$\Tcal_i$ as an element for a given $i \leq k$, they do contain
$\tilde \Tcal_i = \{\tilde q \in \Tcal_i\}$ as an element for each $i \leq k$.
Define $r_k = r$ and
$r_i = r_{i+1} \restriction N_{i+1}(r)$.
Suppose that
$r_{i+1}$ is in $\Tcal_{i+1}$.
Since $\tilde \Tcal_{i+1}$ and $r_{i}$ are in
$N_{i+1}(r) = N_{i+1}(r_{i+1})$ and
since $N_{i+1}(r) \cap H(\aleph_2)$ is in every club in $\Ecal \cap N_{i+1}(r)$,
it follows by elementarity of $N_{i+1}(r)$ that
the set
$$\{N_{i+1}(q^*) \cap H(\aleph_2): q^* \in \Tcal_{i+1} \mand
r_i = q^* \restriction N_{i+1}(q^*)\}=$$
$$\{N_{i+1}(\tilde q^*): \tilde q^* \in \tilde \Tcal_{i+1} \mand
\tilde r_i = \tilde q^* \restriction N_{i+1}(\tilde q^*)\}$$
is stationary.
Hence $r_i$ is in $\Tcal_i$.
\end{proof}

Notice that $\Tcal$ is in $M$.
$\Tcal$ has a natural tree order associated with it induced by
restriction.
Since no element of $\Tcal_k \cap M$ is compatible with $r$ and
since $r_0$ is in $\Tcal \cap M$,
there is a $q$ in $\Tcal \cap M$ which is maximal in the tree order such that
$q$ is compatible with $r$ but such that none of $q$'s immediate
successors in $\Tcal \cap M$ are compatible with $r$.
Let $l$ denote the height of $q$ in $\Tcal$
and put $\Acal$ to be equal to the set of all $N_{l+1}(q^*)$ such that
$q^*$ is an immediate successor of $q$ in $\Tcal$.
Notice that if $q^*$ is in $\Tcal_{l+1}$ and
$q$ is a restriction of $q^*$, then $q^*$ is in $\Tcal$.
Hence we have arranged that
$\{N \cap H(\aleph_2):N \in \Acal\}$ is stationary.
For each $N$ in $\Acal$, select a fixed $q^*$ which is an immediate
successor of $q$ in $\Tcal$
such that $N_{l+1}(q^*) = N$ and put
$$Y(N) = X_{q^*} \setminus X_q.$$
\begin{claim}
For all $N$ in $\Acal \cap M$
the pair $\big(X_r \cup Y(N),\{N\}\big)$ is not a condition in $\partial (K)$.
\end{claim}

\begin{proof}
Let $N$ be in $\Acal \cap M$ and fix an immediate successor
$q^*$ of $q$ in $\Tcal$ such that $N_{l+1}(q^*) = N$ and
$Y(N) = X_{q^*} \setminus X_q$.
Observe that $$(X_r \cup X_{q^*},\Ncal_{q^*} \cup \Ncal_r)$$
is not a condition in $\partial(K)$ but that
$$(X_r \cup X_q,\Ncal_q \cup \Ncal_r)$$ is a condition.
Furthermore,
$(X_r \cup X_{q^*},\Ncal_{q^*} \cup \Ncal_r)$ fails to
be a condition only because it violates item \ref{stationary_positives} in the
definition of $\partial(K)$.
Observe that $\Ncal_{q^*} \setminus \Ncal_q = \{N\}$.
If $N'$ is an element of $\Ncal_r \cup \Ncal_q$, then
the sets of restrictions
$$\{t \restriction (N' \cap \omega_1):t \in X_r \cup X_{q^*}\},$$
$$\{t \restriction (N' \cap \omega_1):t \in X_r \cup X_{q}\}$$
are equal by definitions of $\Tcal_l$ and $q$ and by
our initial assumptions about the closure of
$X_q$ for $q$ in $\Dcal$ under taking certain restrictions.
Since $(X_r \cup X_q, \Ncal_r \cup \Ncal_q)$ is a condition,
such an $N'$ cannot witness the failure of \ref{stationary_positives}.
Therefore it must be the case that the reason
$(X_r \cup X_{q^*},\Ncal_{q^*} \cup \Ncal_r)$ is not in
$\partial (K)$ is that $N$ witnesses a failure of
item \ref{stationary_positives}.
Now, the elements of $X_{q^*}$ which have height at least
$N \cap \omega_1$ are exactly those in $Y(N) = X_{q^*} \setminus X_q$.
This finishes the claim.
\end{proof}
Notice that by elementarity of $M$, $Y \restriction \Acal$
can be chosen to be in $M$.
Now $M$ models ``There is a stationary set of countable elementary
submodels $M_*$ of $H({2^{2^{\aleph_1}}}^+)$ such that for
some $Y(M_*)$ with
$\big(Y(M_*),\{M_* \cap H({2^{\aleph_1}}^+)\}\big)$ in $\partial (K)$
we have that for every $N$ in $\Acal \cap M_*$ the pair
$\big(Y(N) \cup Y(M_*),\{N\}\big)$ is not a condition in $\partial (K)$.''
By elementarity of $M$, we are finished.
\end{proof}

The following definition will be useful.
\begin{defn}
A function $h$ is a \emph{level map} if its domain is a subset of
$\omega_1$ and $h(\delta)$ is a finite subset of the
$\delta\Th$ level of $T$ whenever it is defined.
\end{defn}
The next proposition is useful and follows easily from the fact
that all levels of $T$ are countable.
\begin{prop}\label{level_map_prop}
If $N \cap H(\aleph_2)$ is in $E_0$,
$\delta = N \cap \omega_1$, and $X$ is a finite subset
of the $\delta\Th$ level of $T$, then
there is a level map $h$ in $N$ such that $h(\delta) = X$. 
\end{prop}

The next lemma will represent the only use of $\MRP$ in the proof.

\begin{lem}\label{0-1_law}
($\MRP$)
Suppose that $M$ is a countable elementary submodel of $H({2^{2^{\aleph_1}}}^+)$ which
contains $T$ and $K$ as members.
If $X$ is a finite subset of $T$,
then there is an $E$ in $\Ecal \cap M$ such that either every
element
of $E \cap M$ rejects $X$ or no element of $E \cap M$ rejects $X$.
\end{lem}

\begin{remark}
Notice that the latter conclusion is just a reformulation of
the statement that 
$\big(X,\{M \cap H({2^{\aleph_1}}^+)\}\big)$ is a condition in $\partial(K)$.
\end{remark}

\begin{proof}
Let $\delta = M \cap \omega_1$.
Without loss of generality, we may assume that $X = X \restriction \delta$.
Applying Proposition \ref{level_map_prop},
select a level map $g$ in $M$
such that $g(\delta) = X$.
If $N$ is a countable elementary submodel of
$H({2^{\aleph_1}}^+)$ with $T$ and $K$ as members,
define $\Sigma(N)$ as follows.
If the set of all $P$ in $E_0$ which reject $g(N \cap \omega_1)$
is $N$-stationary, then put $\Sigma(N)$ to be
equal to this set unioned with the complement of $E_0$.
If $\Sigma(N)$ is defined in this way, it will be said to be defined
\emph{non-trivially}.
Otherwise put $\Sigma(N)$ to be the interval
$[\emptyset,N \cap H(\aleph_2)]$.

Observe that $\Sigma$ is an open stationary set mapping which is moreover
an element of $M$.
Applying $\MRP$ and the elementarity of $M$, it is possible to find
a reflecting sequence $\<N_\xi:\xi < \omega_1\>$ for $\Sigma$ which is an
element of $M$.
Let $E$ be the collection of all $\overline P$ in $E_0$ which contain
\begin{enumerate}

\item the sequence $\<N_\xi \cap H(\aleph_2):\xi < \omega_1\>$ and

\item some $\delta_0 < N \cap \omega_1$ such
that $N_\xi \cap H(\aleph_2)$ is in $\Sigma(N_\delta)$
whenever $\xi$ is in $(\delta_0,\delta)$.

\end{enumerate}
Notice that $E$ is in $M \cap \Ecal$.

To finish the proof, suppose that the set of all $P$ in $M \cap E_0$
which reject $X$ is $M$-stationary (i.e. the second conclusion does
not hold).
\begin{claim}
$\Sigma(N_\delta)$ is defined non-trivially.
\end{claim}

\begin{proof}
Suppose that $E' \subseteq E_0$ is a club in $N_\delta$.
Since the reflecting sequence is continuous,
$N_\delta$ is a subset of $M$ and therefore $E'$ is also in $M$.
Let $P$ be an element of $E_0$ with $P \cap \omega_1 = \nu < \delta$.
By assumption, there is a $P$ in $E' \cap M$ such that
$P$ rejects $X$.
Applying elementarily of $N_\delta$, Proposition \ref{rejection_prop}, and the
fact that $X \restriction \nu$ is in $N_\delta$, it is possible
to find such a $P$ in $E' \cap N_\delta$ which rejects
$X \restriction \nu$ --- and hence $X$.
It follows that $\Sigma(N_\delta)$ is defined non-trivially.
\end{proof}

Now suppose that $\overline P$ is in $E \cap M$.
We are finished once we see that $\overline P$ rejects $X$.
Let $\nu = \overline P \cap \omega_1$.
Since $\delta_0 < \nu < \delta$, $P_\nu = N_\nu \cap H(\aleph_2)$ is in $\Sigma(N_\delta)$.
So $P_\nu$ rejects $X$ or --- equivalently ---
$K(X)$ is in $\Ical_{P_\nu}(T)$.
Observe that $P_\nu \cap \omega_1 = \overline P \cap \omega_1$ and
$P_\nu \subseteq \overline P$ by continuity of the reflecting sequence.
Hence $\Ical_P(T) \subseteq \Ical_{\overline P}(T)$.
It follows that $\overline P$ rejects $X$.
\end{proof}

The next lemma finishes the proof of the main theorem.

\begin{lem} \label{MRP_MA_lem} ($\MRP + \MA_{\aleph_1}$)
There are no $\Acal$, $\Bcal$, and $Y$
which satisfy the conclusion of Lemma \ref{A_B:lem}.
In particular, $\partial (K)$ is canonically proper.
\end{lem}

\begin{proof}
We will assume that there are such $\Acal$, $\Bcal$, and $Y$ and derive a
contradiction by violating Lemma \ref{rel_c.c.c.}.
Without loss of generality we may
suppose that elements of $\Bcal$ contain $\Acal$ as a member.
By modifying $Y$ we may assume that
all elements of $Y(M)$ have height $M\cap \omega_1$ whenever
$M$ is in $\Acal \cup \Bcal$.
Further, we may assume that $Y(M)$ has the same fixed $n$ size for all
$M$ in $\Bcal$ and that there is a $\zeta_0$ and $E_* \in \Ecal$ such that:
\begin{enumerate}

\item If $M$ is in $\Bcal$, then
$\diff(Y(M)) \subseteq \zeta_0$.

\item If $M,M'$ are in $\Bcal$, then
$Y(M) \restriction \zeta_0 = Y(M') \restriction \zeta_0$.

\item If $M$ is in $\Bcal$, then $E_*$ is in $M$ and no element of
$E_*$ rejects $Y(M)$.

\end{enumerate}
This is achieved by the pressing down lemma
and the proof of Claim \ref{step_up}.\footnote{To get the last
item, find an $E_1$ in $\Ecal$
such that no element of $E_1 \cap M$ rejects
$Y(M)$ for stationary many
$M$ in $\Bcal$, put $E_*$ to be the elements $P$ of $E_1$ which
are equal to the union of their intersection with $E_1$.}
Let $\Fcal$ be the collection of all finite $X \subseteq T$ such that
all elements of $X$ have the same height $\gamma > \zeta_0$ and the set
$$\{M \in \Bcal : Y(M) \restriction \gamma = X\}$$
is stationary.
Notice that, for a fixed $\gamma$, we can define
$\Bcal$ to be a union over the finite subsets $X$ of $T_\gamma$ of the
collection $$\Bcal[X] = \{M \in \Bcal:Y(M) \restriction \delta = X\}$$
and hence at least one such $\Bcal[X]$ must be stationary.
Consequently $\Fcal$ must be uncountable.
Also, no element of $E_*$ rejects any element of $\Fcal$.
Now define $\Qcal$ to be the collection of all
finite $F \subseteq \Fcal$ such that if
$X\ne X'$ are in $F$,
then the heights of elements of $X$ and $X'$ are different and
there is a $j < n$ such that $X(j) \meet X'(j)$ is not in $K$.

\begin{claim} ($\MRP$)
$\Qcal$ satisfies the countable chain condition.
\end{claim}

\begin{proof}
Suppose that $\<F_\xi : \xi < \omega_1\>$ is a sequence of
distinct elements of $\Qcal$.
We will show that $\{F_\xi:\xi < \omega_1\}$ is not an antichain
in $\Qcal$. 
By a $\Delta$-system argument, we may assume that
the sequence consists of disjoint sets of the same cardinality $m$.
Let $F_\xi(i)$ denote the $i\Th$-least element of $F_\xi$ in the order
induced by $T$'s height function.
If $j < n$, let $F_\xi(i,j)$ denote the $j\Th$ element of
$F_\xi(i)$ in the lexicographical order on $F_\xi(i)$
(i.e. $F_\xi(i,j) = F_\xi(i)(j)$).

Let $N$ be an element of $\Acal$ which contains $\zeta_0$ and
$\<F_\xi:\xi < \omega_1\>$ as members.
Put $\delta = N \cap \omega_1$
and fix a $\beta$ in $\omega_1 \setminus N$.
Let $E$ be a club in $N$ such that
$Y(N)$ is not rejected by any element of $E \cap N$.

For each $i < m$, pick an $M_i$ in $\Bcal$ such that $N \in M_i$ and
$F_\beta(i)$ is a restriction of $Y(M_i)$.
Applying Lemma \ref{0-1_law} for each $i < m$ and intersecting clubs,\footnote{This is
the only place where Lemma \ref{0-1_law} and hence
$\MRP$ is applied.}
it is possible to find a $P$ in $E \cap N$ such that
$$I = \bigcup_{i < m} K\big(Y(N) \cup Y(M_i)\big) \in \Ical_P(T).$$
Put $\nu = P \cap \omega_1$.
Pick a $\zeta < \nu$ such that
$\diff(\cup F_\beta \restriction \nu)$ is contained in $\zeta$
and if $u,v$ are distinct elements of
$\cup F_\beta \restriction \nu$, then $u \restriction \zeta$ and
$v \restriction \zeta$ are distinct.

\begin{subclaim}
There is a sequence $\<\alpha_\xi:\xi < \omega_1\>$ in $P$ such that
for each $\xi < \omega_1$ we have the following conditions:
\begin{enumerate}

\item
$\xi \leq \alpha_\xi$,

\item $\cup F_{\alpha_\xi} \restriction \zeta = \cup F_\beta \restriction \zeta$,

\item $\diff(\cup F_{\alpha_\xi} \restriction \xi) = \diff(\cup F_\beta \restriction \nu)$, and

\item \label{nu_correct}
$\cup F_{\alpha_\nu} \restriction \nu = \cup F_\beta \restriction \nu$.
\end{enumerate}
\end{subclaim}

\begin{proof}
The only part which is non-trivial is to get the sequence to be a
member of $P$ and to satisfy item \ref{nu_correct}.
By Proposition \ref{level_map_prop}, there is a level map $g$ in $P$
such that $g(\nu) = \cup F_\beta \restriction \nu$.
Now working in $P$, we can define $\alpha_\xi$ to be
an ordinal such that $\cup F_{\alpha_\xi} \restriction \xi = g(\xi)$ if
$g(\xi)$ is defined, is a restriction of this form, and satisfies
$\diff(g(\xi)) = \diff(\cup F_\beta \restriction \nu)$ and
$g(\xi) \restriction \zeta = \cup F_\beta \restriction \zeta$.
If $\alpha_\xi$ is left undefined, then simply select a $\alpha_\xi$ with
the necessary properties.
Notice that $\alpha_\nu$ is defined using $g$.
\end{proof}

\begin{subclaim}
There is an uncountable $\Xi \subseteq \omega_1$ in $P$
such that if
$$Z = \{F_{\alpha_\xi}(0,0) \restriction \xi:\xi \in \Xi\}$$ and
$t = F_\beta(0,0) \restriction \nu$, then
$t$ is in the downwards closure of $Z$ and
$\Delta(Z,t)$ is disjoint from 
$$I = \bigcup_{i < m} K\big(Y(N) \cup Y(M_i)\big).$$
\end{subclaim}
\begin{proof}
By Proposition \ref{I_P(T)_refine}
there is a $\Xi_0 \subseteq \omega_1$ such that
for some $t_0$ in $T$ of height $\nu$ in the downward closure of
$Z_0 = \{F_{\alpha_\xi}(0,0) \restriction \xi :\xi \in\Xi_0\}$
the set $\Delta(Z_0,t_0)$ is disjoint from
$I$.
Let $Z_1$ be all elements $s$ in $T$ obtained 
from some $F_{\alpha_\xi}(0,0)\restriction \xi$ by
changing its values on the set
$$\diff\big(F_{\alpha_\nu}(0,0) \restriction \nu,t_0\big) \cap \xi.$$
Let $\Xi$ be the collection of all $\xi$ such that
$F_{\alpha_\xi}(0,0)\restriction \xi$ is an initial
part of some element of $Z_1$.
Notice that $\Xi$ is in $P$ and is uncountable
since it contains $\nu$.
Furthermore, if
$t = F_{\alpha_\nu}(0,0) \restriction \nu = F_{\beta}(0,0) \restriction \nu$,
then $\Delta(Z,t)$ is contained in
$$ \Delta(Z_1,t) = \Delta(Z_0,t_0)$$
and hence is disjoint from $I$.
\end{proof}
The key observation --- and why the main theorem goes through
--- is the following.
Since $P$ is in $E \cap N$, it does not reject $Y(N)$ and
therefore it is the case that
there is a $\xi$ in $\Xi \cap P$
such that for all $j < |Y(N)|$ the restriction
$$Y(N)(j) \restriction \Delta\big(F_{\alpha_\xi}(0,0),t\big)$$
is in $K$ where $t = F_\beta(0,0) \restriction \nu$.
By the choice of $\Xi$ this means that for all $i < m$ there is a
$j < n$ such that
$$F_\beta(i,j) \restriction \Delta\big(F_{\alpha_\xi}(0,0),t\big) =
F_\beta(i,j) \meet F_{\alpha_\xi}(i,j)$$
is not in $K$.
Let $\alpha = \alpha_\xi$.

Now we claim that $F_{\alpha} \cup F_\beta$ is in $\Qcal$.
To see this, suppose that $i,i' < m$.
If $F_\beta(i) \restriction \nu \ne F_\beta(i') \restriction \nu$, then
pick a $j < n$ such that
$F_\alpha(i,j) \meet F_\alpha(i',j)$ is not in $K$.
Since
$$\Delta\big(F_\alpha(i',j),F_\beta(i',j)\big) \geq \zeta >
\Delta\big(F_\alpha(i,j),F_\alpha(i',j)\big)$$
it must be the case that
$$\Delta\big(F_\alpha(i,j),F_\beta(i',j)\big)=
\Delta\big(F_\alpha(i,j),F_\alpha(i',j)\big)$$
and so
$$F_\alpha(i,j) \meet F_\beta(i',j) = F_\alpha(i,j) \meet F_\alpha(i',j)$$
is not in $K$.

If $F_\beta(i) \restriction \nu = F_\beta(i') \restriction \nu$,
then we have that for all $j < n$ that
$$
\Delta\big(F_\alpha(i,j),F_\beta(i',j)\big) =
\Delta\big(F_\alpha(i,j),F_\beta(i,j)\big) =
\Delta\big(F_\alpha(0,0),F_\beta(0,0)\big).$$
By arrangement there is a $j$ such that
$$F_\alpha(i,j) \meet F_\beta(i,j) = F_\beta(i,j) \restriction
\Delta\big(F_\alpha(0,0),F_\beta(0,0)\big)$$
is not in $K$.
Hence for all $i,i' < m$ there is a $j < n$ such that
$$F_\alpha(i,j) \meet F_\beta(i',j)$$ is not in $K$ and therefore
we have that
$F_\alpha \cup F_\beta$ is in $\Qcal$.
\end{proof}

Applying $\MA_{\aleph_1}$ to
the forcing $\Qcal$ it is possible
to find an uncountable $\Fcal_0 \subseteq \Fcal$ such that
whenever $X\ne X'$ are in $\Fcal_0$,
there is a $j < n$ such that $X(j) \meet X'(j)$ is not in $K$.
This contradicts Lemma \ref{rel_c.c.c.} since no element of $E_*$
rejected by any element of $\Fcal$.
\end{proof}

\section{Closing remarks}

The conventional wisdom had been that if it were possible to prove
the consistent existence of a five element basis for the
uncountable linear orders, then such a basis
would follow from $\BPFA$.
The use of $\MRP$ in the argument above
is restricted to proving Lemma \ref{0-1_law}.
Working from a stronger assumption,\footnote{
Questions for the reader:  Why do we need the stronger assumption?
What allows us to use $\MRP$ in the proof of Lemma \ref{0-1_law}?}
the following abstract form of the lemma
can be deduced.
\begin{01law}
($\SMRP$\footnote{
$\SMRP$ is the \emph{Strong Mapping Reflection Principle} obtained
by replacing ``club'' in the statement of $\MRP$ with
``projective stationary'' (see \cite{proj_stat}).
This axiom follows from Martin's Maximum
via the same proof that $\MRP$ follows from $\PFA$ (see \cite{MRP}).})
Suppose that $\Sigma$ is an open set mapping defined on a club
and that $\Sigma$ has the following properties:
\begin{enumerate}

\item If $N$ is in the domain of $\Sigma$, then
$\Sigma(N)$ is closed under end extensions.\footnote{
Here we define $\overline N$ \emph{end extends} $N$ as meaning that
$N \cap \omega_1 = \overline N \cap \omega_1$ and $N \subseteq \overline N$.}

\item If $N$ and $\overline N$ are in the domain of $\Sigma$ and
$\overline N$ is an end extension of $N$, then
$\Sigma(N) = \Sigma(\overline N) \cap N$.

\end{enumerate}
Then for a closed unbounded set of $N$ in the domain of $\Sigma$,
there is a club $E \subseteq [X_\Sigma]^{\aleph_0}$ in $N$
such that $E \cap N$ is either contained in or disjoint from
$\Sigma(N)$.
\end{01law}
\noindent
It seems quite possible that this \emph{0-1 law}
will be useful in analyzing related problems such as
Fremlin's problem on perfectly normal compacta
(see \cite{open_problems_topology:Gruenhage}, \cite{basis_problems}).

$\MRP$ has considerable consistency strength
\cite{MRP}, while
$\BPFA$ can be forced if there is a reflecting cardinal
\cite{BPFA}.
The following is left open.
\begin{question}
Does $\BPFA$ imply Shelah's conjecture?
\end{question}

Recently K\"onig, Larson, Veli\v{c}kovi\'c, and I have shown that
a certain saturation property of Aronszajn trees taken
together with $\BPFA$ implies Shelah's conjecture.
This saturation property can be forced if there is a Mahlo cardinal.
This considerably reduces the upper bound on the consistency strength
of Shelah's conjecture to that of a reflecting Mahlo cardinal.
It is possible, however, that Shelah's conjecture cannot follow
from $\BPFA$ simply on grounds of its consistency strength.
\begin{question}
Does Shelah's conjecture imply that $\aleph_2$
is either Mahlo or reflecting in $L$?
\end{question}
Reflecting cardinals are weaker in consistency strength than Mahlo cardinals;
if the proper class ordinal is Mahlo, then there is a proper
class of reflecting cardinals.
Still, the least Mahlo cardinal is not reflecting.
It should be remarked though that
Shelah's conjecture is not known to have any
large cardinal strength.


\begin{thebibliography}{10}

\bibitem{club_isomorphic}
U.~Abraham and S.~Shelah.
\newblock Isomorphism types of {A}ronszajn trees.
\newblock {\em Israel J. Math.}, 50(1-2):75--113, 1985.

\bibitem{probabilistic_method}
Noga Alon and Joel~H. Spencer.
\newblock {\em The probabilistic method}.
\newblock Wiley-Interscience Series in Discrete Mathematics and Optimization.
  Wiley-Interscience [John Wiley \& Sons], New York, second edition, 2000.

\bibitem{gen_abs_FA}
Joan Bagaria.
\newblock Generic absoluteness and forcing axioms.
\newblock In {\em Models, algebras, and proofs (Bogot\'a, 1995)}, pages 1--12.
  Dekker, New York, 1999.

\bibitem{reals_isomorphic}
James~E. Baumgartner.
\newblock All {$\aleph \sb{1}$}-dense sets of reals can be isomorphic.
\newblock {\em Fund. Math.}, 79(2):101--106, 1973.

\bibitem{set_theory:Bekkali}
M.~Bekkali.
\newblock {\em Topics in set theory}.
\newblock Springer-Verlag, Berlin, 1991.
\newblock {L}ebesgue measurability, large cardinals, forcing axioms,
  $\rho$-functions, Notes on lectures by {S}tevo {T}odor\v cevi\'c.

\bibitem{proj_stat}
Qi~Feng and Thomas Jech.
\newblock Projective stationary sets and a strong reflection principle.
\newblock {\em J. London Math. Soc. (2)}, 58(2):271--283, 1998.
arXiv:math.LO/9409202.

\bibitem{BPFA}
Martin Goldstern and Saharon Shelah.
\newblock The {B}ounded {P}roper {F}orcing {A}xiom.
\newblock {\em J. Symbolic Logic}, 60(1):58--73, 1995.
GoSh:507. arXiv:math.LO/9501222.

\bibitem{open_problems_topology:Gruenhage}
Gary Gruenhage.
\newblock Perfectly normal compacta, cosmic spaces, and some partition
  problems.
\newblock In {\em Open problems in topology}, pages 85--95. North-Holland,
  Amsterdam, 1990.

\bibitem{multiple_forcing}
T.~Jech.
\newblock {\em Multiple forcing}, volume~88 of {\em Cambridge Tracts in
  Mathematics}.
\newblock Cambridge University Press, Cambridge, 1986.

\bibitem{set_theory:Jech}
Thomas Jech.
\newblock {\em Set theory}.
\newblock Perspectives in Mathematical Logic. Springer-Verlag, Berlin, second
  edition, 1997.

\bibitem{set_theory:Kunen}
Kenneth Kunen.
\newblock {\em An introduction to independence proofs}, volume 102 of {\em
  Studies in Logic and the Foundations of Mathematics}.
\newblock North-Holland, 1983.

\bibitem{MRP}
Justin~Tatch Moore.
\newblock Set mapping reflection.
\newblock submitted to JML in Nov. 2003.

\bibitem{Countryman:Shelah}
Saharon Shelah.
\newblock Decomposing uncountable squares to countably many chains.
\newblock {\em J. Combinatorial Theory Ser. A}, 21(1):110--114, 1976.

\bibitem{proper_forcing}
Saharon Shelah.
\newblock {\em Proper and improper forcing}.
\newblock Springer-Verlag, Berlin, second edition, 1998.

\bibitem{problems:Shelah}
Saharon Shelah.
\newblock On what {I} do not understand (and have something to say). {I}.
\newblock {\em Fund. Math.}, 166(1-2):1--82, 2000.
\newblock Saharon Shelah's anniversary issue. Sh:666.
arXiv:math.LO/9906113.

\bibitem{many_iso_types}
W.~Sierpi\'{n}ski.
\newblock Sur un probl\`eme concernant les types de dimensions.
\newblock {\em Fundamenta Mathematicae}, 19:65--71, 1932.

\bibitem{trees:Todorcevic}
Stevo Todor{\v{c}}evi{\'c}.
\newblock Trees and linearly ordered sets.
\newblock In {\em Handbook of set-theoretic topology}, pages 235--293.
  North-Holland, Amsterdam, 1984.

\bibitem{basis_problems}
Stevo Todor{\v{c}}evi{\'c}.
\newblock Basis problems in combinatorial set theory.
\newblock In {\em Proceedings of the International Congress of Mathematicians,
  Vol. II (Berlin, 1998)}, number Extra Vol. II, pages 43--52, 1998.

\bibitem{l-map}
Stevo Todor\v{c}evi\'{c}.
\newblock Lipszhitz maps on trees.
\newblock report 2000/01 number 13, Institut Mittag-Leffler.

\bibitem{partitioning_ordinals}
Stevo Todor\v{c}evi\'{c}.
\newblock Partitioning pairs of countable ordinals.
\newblock {\em Acta Math.}, 159(3--4):261--294, 1987.

\bibitem{partition_problems}
Stevo Todor\v{c}evi\'{c}.
\newblock {\em Partition Problems In Topology}.
\newblock Amer. Math. Soc., 1989.

\bibitem{class_trans}
Stevo Todor\v{c}evi\'{c}.
\newblock A classification of transitive relations on {$\omega\sb 1$}.
\newblock {\em Proc. London Math. Soc. (3)}, 73(3):501--533, 1996.

\bibitem{cseq}
Stevo Todor\v{c}evi\'{c}.
\newblock Coherent sequences.
\newblock In {\em Handbook of Set Theory}. North-Holland, (in preparation).

\bibitem{OCA_aut}
Boban Veli\v{c}kovi\'{c}.
\newblock {O}{C}{A} and automorphisms of ${\Pcal}(\omega)/\fin$.
\newblock {\em Topology Appl.}, 49(1):1--13, 1993.

\end{thebibliography}
\end{document}